\documentclass[a4paper,11pt, final]{article}

\usepackage{graphicx}
\usepackage{amssymb,amsmath,amsthm,mathrsfs,xfrac}
\usepackage[shortlabels]{enumitem}
\usepackage{a4wide}
\usepackage{hyperref}
\usepackage{array}
\usepackage{colonequals}

\usepackage{tikz}

\numberwithin{equation}{section}
\allowdisplaybreaks[4]

\newtheorem{proposition}{Proposition}[section]
\newtheorem{theorem}[proposition]{Theorem}

\newtheorem{lemma}[proposition]{Lemma}
\newtheorem{definition}[proposition]{Definition}
\newtheorem{corollary}[proposition]{Corollary}
\newtheorem{remark}[proposition]{Remark}

\newenvironment{proofof}[1]{\smallskip\noindent{\textbf{Proof~of~#1.}}%
  \hspace{1pt}}{\hspace{-5pt}{\nobreak\quad\nobreak\hfill\nobreak%
    $\square$\vspace{2pt}\par}\smallskip\goodbreak}

\newcommand{\pint}[1]{\mathaccent23{#1}}
\newcommand{\C}[1]{\mathbf{C}^{#1}}

\newcommand{\BV}{\mathbf{BV}}
\newcommand{\SBV}{\mathbf{SBV}}

\newcommand{\modulo}[1]{{\left|#1\right|}}

\newcommand{\reali}{{\mathbb{R}}}

\newcommand{\tv}{\mathinner{\rm TV}}

\newcommand{\Lip}{\mathinner\mathbf{Lip}}
\renewcommand{\L}[1]{{\mathbf{L}^#1}}

\renewcommand{\epsilon}{\varepsilon}
\renewcommand{\phi}{\varphi}
\renewcommand{\theta}{\vartheta}

\renewcommand{\d}[1]{\mathinner{\mathrm{d}{#1}}}

\newcommand{\essinf}{\mathinner{\mathop{{\rm ess\,inf}}}}
\newcommand{\esssup}{\mathinner{\mathop{{\rm ess\,sup}}}}

\newcommand{\claimend}{\hfill$\checkmark$
  \smallskip\goodbreak}

\renewcommand{\SS}{S^{{\scriptscriptstyle HJ}}}
\newcommand{\II}{I^{\scriptscriptstyle HJ}}
\renewcommand{\ss}{S^{\scriptscriptstyle CL}}
\newcommand{\ii}{I^{\scriptscriptstyle CL}}

\newcommand\rest[2]{{#1}_{\mkern 1mu \vrule height 2ex \mkern 2mu {#2}}}

\DeclareFontFamily{U}{mathx}{}
\DeclareFontShape{U}{mathx}{m}{n}{<-> mathx10}{}
\DeclareSymbolFont{mathx}{U}{mathx}{m}{n}
\DeclareMathAccent{\widecheck}{0}{mathx}{"71}

\title{Localized Inverse Design in\\
  Conservation Laws and Hamilton-Jacobi Equations}

\author{Rinaldo M.~Colombo$^1$ \and Vincent Perrollaz$^2$}

\begin{document}

\maketitle

\footnotetext[1]{INdAM Unit, University of Brescia,
  Italy. \texttt{rinaldo.colombo@unibs.it}}

\footnotetext[2]{Institut Denis Poisson, Universit\'e de Tours, CNRS
  UMR 7013, Universit\'e d'Orl\'eans, France\\
  \texttt{vincent.perrollaz@univ-tours.fr}}

\begin{abstract}

  \noindent Consider the \emph{inverse design} problem for a scalar
  conservation law, i.e., the problem of finding initial data evolving
  into a given profile at a given time. The solution we present below
  takes into account \emph{localizations} both in the final interval
  where the profile is assigned and in the initial interval where the
  datum is sought, as well as additional \emph{a priori}
  \emph{constraints} on the datum's range provided by the model. These
  results are motivated and can be applied to data assimilation
  procedures in traffic modeling and accidents localization.

  \medskip

  \noindent\textit{2000~Mathematics Subject Classification:} 35L65,
  90B20, 35R30.

  \medskip


\end{abstract}


\section{Introduction}
\label{sec:Intro}

We deal with the -- constrained and localized -- inverse design
related to the Cauchy problems for the scalar one dimensional
conservation law and for the Hamilton-Jacobi equation. These equations
read
\begin{equation}
  \label{eq:28}
  \left\{
    \begin{array}{l}
      \partial_t u + \partial_x f (u) = 0
      \\
      u (0,x) = u_o (x)
    \end{array}
  \right.
  \qquad \mbox{ and } \qquad
  \left\{
    \begin{array}{l}
      \partial_t U + f (\partial_x U) = 0
      \\
      U (0,x) = U_o (x)
    \end{array}
  \right.
\end{equation}
where $t \in \reali_+$ is time and $x \in \reali$ is the space
coordinate. In~\eqref{eq:28}, on the left, $u = u (t,x)$ is the
unknown density of a conserved variable and $f$ is the flux, while on
the right $U = U (t,x)$ is the unknown value function and $f$ is the
Hamiltonian. The Cauchy problems~\eqref{eq:28} generate the semigroups
\begin{equation}
  \label{eq:19}
  \ss \colon
  \reali^+ \times \L\infty (\reali; \reali)
  \to
  \L\infty (\reali; \reali)
  \quad \mbox{ and } \quad
  \SS \colon
  \reali^+ \times \Lip (\reali; \reali)
  \to
  \Lip (\reali; \reali)
\end{equation}
in the sense that the orbits $t \to \ss_t u_o$ and $t \to \SS_t U_o$
solve~\eqref{eq:28}, respectively in the entropy or viscosity sense.
For a given real interval $J$ -- the \emph{constraint} -- we provide a
full characterization of the \emph{constrained inverse designs} for
the two equations in~\eqref{eq:28}, namely
\begin{equation}
  \label{eq:5}
  \begin{array}{@{}r@{\,}c@{\,}l}
    \ii_T(u_T;J)
    & \colonequals
    & \left\{
      u_o \in \L\infty (\reali; \reali) \colon
      u_o (\reali) \subseteq J
      \mbox{ and }
      \ss_T u_o = u_T
      \right\}
    \\ \II_T(U_T;J)
    & \colonequals
    & \left\{
      U_o \in \Lip (\reali; \reali) \colon
      U_o' (\reali) \subseteq J
      \mbox{ and }
      \SS_T U_o = U_T
      \right\}
  \end{array}
\end{equation}
for given functions $u_T \in \L\infty (\reali; J)$,
$U_T \in \Lip (\reali; \reali)$ with $U_T' (\reali) \subseteq J$ and
for a fixed $T > 0$.

Then, we extend this characterization to comprehend cases where the
profile $u_T$ is known only on a given real interval, say
$K_T$:
\begin{displaymath}
  \rest{\ii_T (u_T;J)}{K_o}
  \colonequals
  \left\{
    \tilde u_o \in \L\infty (K_o;J) \colon
    \exists \, u_o \in \L\infty (\reali; J)
    \mbox{ with }
    \begin{array}{l}
      \rest{u_o}{K_o}
      =
      \tilde u_o
      \mbox{ and }
      \\
      \rest{\ss_T u_o}{K_T}
      =
      u_T
    \end{array}
  \right\}.
\end{displaymath}
and the meaningful interval $K_o$ is singled out below by means of
$K_T$ and $u_T$.

Our motivation is a typical situation in traffic management: there,
traffic flow measurements are available at a given location, say
$x=L$, during a given time interval. Out of these data, one seeks to
reconstruct the flow along road segments before and after the
measuring site.  Thus, one is lead to solve an inverse design problem,
with constraints on the unknown function (traffic density varies in
the fixed bounded interval $J$) and localized in space time, say at
$\{L\} \times K_T$. Refer to Section~\ref{sec:traffic} for the
detailed discussion.

Data assimilation and flow reconstruction are a common problem in
various disciplines: the monograph~\cite{kalnay2002atmospheric} deals
with the case of weather forecasts, a field classically related to
these problems. Applications to oil reservoirs are
in~\cite{WOS:000305191600013} while the special
issue~\cite{Suzuki2015} is devoted to general fluid dynamics
and~\cite{MR4344886, MR1804497} provide a more analytical approach
devoted to Navier--Stokes equation. We refer to~\cite{MR3363508} for
further examples. Deeply related to the present result
is~\cite{zbMATH06445464}, where kinetic techniques are employed.

The current literature provides several results about the inverse
design for conservation laws. The case of the (inviscid) Burgers'
equation is thoroughly considered in~\cite{MR3643881, MR4437890,
  MR4598825}, the general homogeneous case is solved
in~\cite{ColomboPerrollaz} while the $x$ dependent case is tackled
in~\cite{ColomboPerrollazSylla2}, see also~\cite{MR4503821} for an
alternative approach and~\cite{MR4173218} for results on the Hamilton
- Jacobi equation in several space dimensions. A specific system is
considered in~\cite{coclite:hal-04164795}.

The above applications and, in particular, our motivation are based on
conservation laws. However, we extensively make use of techniques
typical of Hamilton-Jacobi equation and we extensively exploit the
deep connection between the two classes of equations, summarized
in~\eqref{eq:10} below. This allows, in particular, to obtain
Theorem~\ref{thm:flat}, which provides a new independent --- and simpler --- proof of the
characterization of the reachability of a profile also in the case of
data constrained to attain values in $J$. This result improves that
in~\cite{ColomboPerrollazSylla2} relaxing the regularity assumptions
on the target profile.

As the statements and proofs below show, the roles of the constraint
$J$ and of the localization to $K_T$ are entirely different. The
former one appears as somewhat marginal, see
Remark~\ref{rem:label}. This is due to our requiring that the flux $f$
in~\eqref{eq:28} is convex and independent of $x$. Indeed, the example
in~\cite{ColomboPerrollazSylla2, example} shows that, without these
conditions, the range of the solution to a conservation law may well
significantly grow.

When dealing with real measurements, it is necessary to further
localize the available information since, in general, only discrete
samples of data can be collected. Moreover, measurement errors affect
these values. A refinement of the analytical techniques below might
tackle these practical issues.

\medskip

The next section presents our results, first only in the constrained
case, then also under a localization
condition. Section~\ref{sec:traffic} is devoted to the implications of
the general results to our original motivation rooted in vehicular
traffic modeling. Finally, all proofs are deferred to
Section~\ref{sec:Tech}.

\nocite{MR4598825}

\section{Analytic Results}
\label{subs:Notation}

Throughout, concerning conservation laws, we refer to entropy
solutions as classically defined by {Kru\v
  zkov}~\cite[Definition~1]{Kruzhkov1970}, see
also~\cite[Definition~2.1]{ColomboPerrollaz}. In the case of the
Hamilton-Jacobi equation we refer to viscosity
solutions~\cite[Definition~I.1]{CL1983}, see
also~\cite[Definition~2.2]{ColomboPerrollaz}. The connection between
the Cauchy problems in~\eqref{eq:28} is summarized by the following
commuting diagrams:
\begin{equation}
  \label{eq:10}
  \begin{array}{c@{\qquad\qquad\qquad}c}
    \mbox{\cite[Theorem~1.1]{MR1923522}\ }
    &\mbox{\cite[Proposition~2.3]{ColomboPerrollaz}\qquad\qquad\qquad\quad}
    \\
    \begin{array}{c@{}ccc@{}c}
      & U_o
      & \longrightarrow
      & \SS_t U_o
      \\
      \partial_x
      & \big\downarrow
      &
      & \big\downarrow
      & \partial_x
      \\
      & u_o
      & \longrightarrow
      & \ss_t u_o
    \end{array}
      & \begin{array}{c@{}ccc@{}c}
          & U_o
          & \longrightarrow
          & \SS_t U_o
          \\
          \int^x
          & \big\uparrow
          &
          & \big\uparrow
          & \mbox{\cite[Formula~(2.2)]{ColomboPerrollaz}}
          \\
          & u_o
          & \longrightarrow
          & \ss_t u_o
        \end{array}
  \end{array}
\end{equation}
see also~\cite{ColomboPerrollazDepX, MR1923522}.

On the flux/Hamiltonian $f$ we require the following condition:
\begin{enumerate}[label={\bf (f)}]
\item \label{item:1} $f \in \C2 (\reali ; \reali)$ is strongly convex,
  in the sense that $f'' (x) > 0$ for all $x \in \reali$,
  $\lim_{x\to-\infty} f' (x) = -\infty$ and
  $\lim_{x\to+\infty} f' (x) = +\infty$.
\end{enumerate}

Throughout, for the basic results about Hamilton-Jacobi equation we
refer to~\cite[Theorem~4 and Theorem~5, Section~3,
Chapter~3]{Evans}. In particular, the solution to~\eqref{eq:4} can be
written as
\begin{equation}
  \label{eq:17}
  \forall\, T \in \reali_+ \setminus\{0\} \quad
  \forall\, x \in \reali \qquad
  (\SS_T U_o)(x)
  =
  \inf_{\xi \in \reali} \left[
    U_o (\xi) + T \; f^*\left(\frac{x-\xi}{T}\right)
  \right]
\end{equation}
and the Legendre transform $f^*$ of $f$ is recalled in
Definition~\ref{def:Legendre}.

For fixed $u \in \L\infty (\reali;\reali)$ and $T>0$, from the theory
of conservation laws we introduce
\begin{equation}
  \label{eq:15}
  \begin{array}{ccccl}
    \pi_u
    & \colon
    & \reali
    & \to
    & \reali
    \\
    &
    & x
    & \mapsto
    & x - T \, f' \! \left(u (x)\right) \,.
  \end{array}
\end{equation}
As soon as $u$ is the solution to the conservation law~\eqref{eq:28}
at time $T$, $\pi_u (x)$ is the intersection between the axis $t = 0$
and the minimal --- for $u$ left continuous at $x$ --- backward
characteristic from $(T, x)$, see~\cite[\S~10.3 and
\S~11.1]{MR3468916}.

Fix $T>0$. We say that a map $u \in \L\infty (\reali;\reali)$
satisfies Oleinik estimate~\cite{Oleinik1957Originale} at $T$,
whenever the following condition hold:
\begin{enumerate}[label={\bf(O)}]
\item \label{item:2} For a.e.~$x \in \reali$ and
  $\Delta x \in \reali_+ \setminus\{0\}$,
  $\dfrac{f'\left(u (x+\Delta x)\right) - f'\left(u (x)\right)}{\Delta
    x} \leq \dfrac{1}{T}$.
\end{enumerate}
As is well known, $u$ satisfies condition~\ref{item:2} at $T$ if and
only if the map $\pi_u$ is (a.e.) non decreasing. Hence, since
$u (x) = (f')^{-1} \left(\frac{x-\pi_u (x)}{T}\right)$,
by~\ref{item:1} $u$ admits a representative which is left continuous
and thus Condition~\ref{item:2} is satisfied for \emph{every}
$x \in \reali \setminus \{0\}$ and $\Delta x \in \reali_+$. Below,
when~\ref{item:2} applies, by $u$ or $u_T$ we understand their left
continuous representatives.

\subsection{Constrained Inverse Design}
\label{subsec:inverse-design}

For a non empty closed real interval $J$ acting as \emph{constraint}
we now head towards detailed descriptions of the sets~\eqref{eq:5} of
profiles $u_T$ and $U_T$ that can be attained as solutions
to~\eqref{eq:28} for suitable initial data $u_o$ and $U_o$.  The case
$J = \reali$ is not excluded.  We stress that the results below depend
on and take care of the constraint $J$. By~\eqref{eq:10}, the
descriptions of the two sets are one consequence of the other.

Proceeding towards a full characterization of $\II_T (U_T; J)$, we
start collecting information on a \emph{particular} element
$U_o^\flat$.

Proofs to statements in this paragraph are deferred to
\S~\ref{subsec:proofs-inverse-design}.

\begin{theorem}
  \label{thm:flat}
  Let $f$ satisfy~\ref{item:1}, $J$ be a non empty closed interval and
  $T$ be positive. Fix $u_T \in \L\infty (\reali; J)$ satisfying
  Condition~\ref{item:2} at $T$. Define, for a fixed
  $\widecheck x \in \reali$, for all $x \in \reali$
  \begin{equation}
    \label{eq:7}
    U_T (x)
    \colonequals
    \int_{\widecheck x}^x u_T (\xi)\d\xi
    \,,\;
    U_o^\flat (x)
    \colonequals
    \sup_{\xi \in \reali} \left[
      U_T (\xi)
      - T \; f^* \! \left(\dfrac{\xi-x}{T}\right)
    \right]
    \,,\;
    u_o^\flat (x) \colonequals \frac{\d~}{\d x} U_o^\flat (x) \,.
  \end{equation}
  Then, $U_o^\flat \in \II_T (U_T; J)$ and
  $u_o^\flat \in \ii_T (u_T; J)$.
\end{theorem}

Above, $f^*$ is, as usual, the Legendre transform of $f$, see
Definition~\ref{def:Legendre}.

\begin{corollary}
  \label{prop:Oleinik}
  Let $f$ satisfy~\ref{item:1}, $J$ be a non empty closed interval and
  $T$ be positive. Fix $u_T \in \L\infty (\reali; J)$. Then, $u_T$
  satisfies Condition~\ref{item:2} at $T$ if and only if
  $\ii_T (u_T; J) \ne \emptyset$.
\end{corollary}

\noindent The proof of Corollary~\ref{prop:Oleinik} is as follows: on
the one hand, Oleinik's result~\cite{Oleinik1957Originale} ensures
that Condition~\ref{item:2} holds for any reachable profile. On the
other hand, Theorem~\ref{thm:flat} shows that Condition~\ref{item:2}
implies that $\ii_T (u_T; J)$ is non empty.

Corollary~\ref{prop:Oleinik} was originally stated
as~\cite[Corollary~3.2]{ColomboPerrollaz} in the case $J=\reali$, and
proved relying on conservation laws techniques, while the proof of
Theorem~\ref{thm:flat} is based on the Hamilton-Jacobi
equation~\eqref{eq:4}.

An explicit construction of $u_o^\flat$ is provided in the following
extension of~\cite[Theorem~3.1]{ColomboPerrollaz} which, besides
comprising the constraint $J$, does not require $\SBV$ regularity.

\begin{proposition}
  \label{prop:reverse}
  Let $f$ satisfy~\ref{item:1} and $T$ be positive. Fix
  $u_T \in \L\infty (\reali; \reali)$. Define $U_T$, $U_o^\flat$ and
  $u_o^\flat$ as in~\eqref{eq:7}. Call $v$ and $V$ the entropy and
  viscosity solutions to
  \begin{equation}
    \label{eq:29}
    \left\{
      \begin{array}{l}
        \partial_t v + \partial_x f (-v) = 0
        \\
        v (0,x) = -u_T (x)
      \end{array}
    \right.
    \quad \mbox{ and } \quad
    \left\{
      \begin{array}{l}
        \partial_t V + f (-\partial_x V) = 0
        \\
        V (0,x) = -U_T (x) \,.
      \end{array}
    \right.
  \end{equation}
  Then, for a.e.~$t \in [0,T]$ and $x \in \reali$,
  \begin{displaymath}
    \begin{array}{rcl}
      v (T,x)
      & =
      & -v_o^\flat (x)\,;
      \\
      V (T,x)
      & =
      & -U_o^\flat (x)\,;
    \end{array}
    \qquad
    v (t,x) = \partial_x V (t,x)
    \quad \mbox{ and } \quad
    \tv (u_T) \geq \tv (u_o^\flat) \,.
  \end{displaymath}
\end{proposition}

When $f$ depends explicitly on $x$, the bending of characteristic lines
significantly complicates a result like
Proposition~\ref{prop:reverse}, see~\cite{ColomboPerrollazSylla2}.

The function $u_o^\flat$ defined in~\eqref{eq:7} has the minimal range
among those profiles that evolve into $u_T$, as proved by the
following result.

\begin{corollary}
  \label{cor:clco}
  Let $f$ satisfy~\ref{item:1} and $T$ be positive. Fix
  $u_T \in \L\infty (\reali; \reali)$ satisfying
  Condition~\ref{item:2} at $T$. Let $u_o^\flat$ be as defined
  in~\eqref{eq:7}. Then,
  \begin{displaymath}
    \tv (u_T) = \tv (u_o^\flat)
    \quad \mbox{ and } \quad
    \overline{\mathop{\rm co}} \; u_T (\reali)
    =
    \overline{\mathop{\rm co}} \; u_o^\flat (\reali) \,.
  \end{displaymath}
\end{corollary}

\noindent Above, for $A \subseteq \reali$,
$\overline{\mathop{\rm co}} \; A$ stands for the closed convex hull of
$A$.

\begin{remark}
  \label{rem:label}
  \rm Corollary~\ref{cor:clco} underlines that the role of $J$ in
  reachability is rather marginal, as soon as it is sufficiently
  large, i.e., as soon as
  $J \supseteq \overline{\mathop{\rm co}} \; u_T (\reali)$.
\end{remark}

\begin{proposition}[{\cite[Theorem~11.4.3]{MR3468916}}]
  \label{thm:two_min_renewed}
  Let $f$ satisfy~\ref{item:1} and $T$ be positive. For all
  $U_o \in \Lip (\reali;\reali)$, define $u_o = U_o'$ and set
  $u_T = \ss_T u_o$, according to~\eqref{eq:19}. Using the
  notation~\eqref{eq:15}, for all $x \in \reali$
  \begin{equation}
    \label{eq:18}
    \SS_T U_o (x)
    =
    U_o\left(\pi_{u_T} (x)\right)
    + T \; f^*\left(\frac{x-\pi_{u_T} (x)}{T}\right) \,.
  \end{equation}
\end{proposition}

Note that, when $f$ also depends explicitly on $x$,
\cite[Theorem~3.1]{ColomboPerrollazSylla2} reformulates
Proposition~\ref{thm:two_min_renewed} in terms of the connection
between generalized characteristics and minima of the integral
functional connected to the Hamilton-Jacobi equation.

\begin{proposition}
  \label{prop:characterization}
  Let $f$ satisfy~\ref{item:1}, $J$ be a non empty closed interval and
  $T$ be positive. Fix $u_T \in \L\infty (\reali; J)$ satisfying
  Condition~\ref{item:2} at $T$. Define $U_T$ and $U_o^\flat$ as
  in~\eqref{eq:7}.  Then, for all $u_o \in \L\infty (\reali; \reali)$,
  using the notation~\eqref{eq:15} and setting, for any fixed
  $\widecheck x \in \reali$,
  \begin{equation}
    \label{eq:11}
    \forall\,t \in \reali \qquad
    U_o (x)
    \colonequals
    \int_{\pi_{u_T} (\widecheck x)}^x u_o (\xi) \d\xi -
    T\, f^*\left(-\frac{\pi_{u_T} (\widecheck x)}{T}\right)
  \end{equation}
  we have the equivalences
  \begin{equation}
    \label{eq:3}
    u_o \in \ii_T (u_T; J)
    \iff
    U_o \in \II_T (U_T; J)
    \iff
    \left\{
      \begin{array}{rl}
        (i)
        & U_o \geq U_o^\flat \,;
        \\
        (ii)
        & U_o = U_o^\flat \mbox{ on } \overline{\pi_{U'_T} (\reali)} \,;
        \\
        (iii)
        & U_o' (\reali) \subseteq J \,.
      \end{array}
    \right.
  \end{equation}
\end{proposition}

The above proposition is strictly related
to~\cite[Theorem~3.3]{ColomboPerrollazSylla2}
and~\cite[Theorem~2.6]{MR4173218}: both these results do not consider
the constraint $J$ but the former one applies to general $x$ dependent
Hamiltonian functions while the latter applies to several space
dimensions.

When the constraint $J$ is compact, the introduction of the maps
$U^\sharp_o$ and $u^\sharp_o$ in the following Proposition allows the
precise description of the inverse design sets provided by
Theorem~\ref{thm:1}.

\begin{proposition}
  \label{prop:sharp}
  Let $f$ satisfy~\ref{item:1}, $J$ be a non empty compact interval
  and $T$ be positive. Fix $u_T \in \L\infty (\reali; J)$
  satisfying~\ref{item:2} at $T$. Define $U_T$ as in~\eqref{eq:7} and
  \begin{equation}
    \label{eq:6}
    \forall\,x \in \reali \qquad
    U_o^\sharp (x)
    \colonequals
    \sup \left\{
      U_o (x) \colon U_o \in \II_T (U_T; J)
    \right\} \,,
    \qquad
    u_o^\sharp (x)  \colonequals \frac{\d~}{\d x} U_o^\sharp (x) \,.
  \end{equation}
  Then, $U_o^\sharp \in \II_T (U_T; J)$ and
  $u_o^\sharp \in \ii_T (u_T; J)$.
\end{proposition}

\noindent Note that while $U_T$ and $U^\sharp_o$ depend on
$\widecheck x$ in~\eqref{eq:7}, the map $u^\sharp_o$ is actually
independent of it.

We are now ready for a further characterization of the set
$\II_T (U_T; J)$.

\begin{theorem}
  \label{thm:1}
  Let $f$ satisfy~\ref{item:1}, $J$ be a non empty compact interval
  and $T$ be positive. Fix $u_T \in \L\infty (\reali; J)$
  satisfying~\ref{item:2} at $T$. Then, with the
  notation~\eqref{eq:7}--\eqref{eq:6},
  \begin{equation}
    \label{eq:14}
    \II_T (U_T; J)
    =
    \left\{
      U_o \in \Lip (\reali; \reali)
      \colon
      \begin{array}{l}
        U_o' (x)
        \in
        J
        \\
        U_o (x) \in [U_o^\flat (x), U_o^\sharp (x)]
      \end{array}
      \mbox{ for a.e. } x \in \reali
    \right\} \,.
  \end{equation}
  In particular, $\II_T (U_T; J)$ is convex and compact with respect
  to the topology of uniform convergence on compact subsets of
  $\reali$. Moreover, setting
  $\widecheck y \colonequals \pi_{u_T} (\widecheck x)$, for any fixed
  $\widecheck x \in \reali$, we have the characterization
  \begin{eqnarray}
    \!\!\!\!\!\!
    \ii_T (u_T;J)
    \label{eq:16}
    & =
    & \left\{
      u_o \in \L\infty (\reali; J)
      \colon
      \textstyle
      \int_{\widecheck y}^y  u_o \d{x}
      \in \left[
      \int_{\widecheck y}^y  u_o^\flat \d{x},
      \int_{\widecheck y}^y  u_o^\sharp \d{x}
      \right]
      \mbox{ for all }y \in \reali
      \right\}
  \end{eqnarray}
  so that $\ii_T (u_T)$ is convex and sequentially compact with
  respect to the weak-$*$ $\L\infty$ topology.
\end{theorem}

\noindent Again, note that the arbitrariness of $\widecheck x$ does
affect $U^\flat_o$ and $U^\sharp_o$ bus has no relevance on
$u^\flat_o$, $u^\sharp_o$.

\subsection{Inverse Design Localized in Space}
\label{subsec:time-localization}

This section is devoted to the localization of the previous results on
two space intervals: the former one, $K_T$, is to be considered at
time $t=T$ and the latter one, $K_o$ at time $t=0$.

Hereafter, we consider only the conservation law in~\eqref{eq:28}, the
case of the Hamilton-Jacobi equation being entirely analogous.

The proofs related to statements in this section, where necessary, are
deferred to \S~\ref{subsec:proofs-time-localization}.

\begin{definition}
  \label{def:reachable}
  Let $J$ be a non trivial closed real interval and $T>0$. Fix a
  second non trivial closed real interval $K_T$. A profile
  $u_T \in \L\infty (K_T; J)$ is \emph{reachable at $t=T$ on $K_T$} if
  there exists a $u_o \in \L\infty (\reali; J)$ such that the
  corresponding solution $u$ to the conservation law in~\eqref{eq:28}
  satisfies $\rest{\ss_Tu_o}{K_T} = u_T$.

  If $K_o$ is another non trivial closed real interval, denote
  \begin{equation}
    \label{eq:30}
    \rest{\ii_T (u_T;J)}{K_o}
    \colonequals
    \left\{
      \tilde u_o \in \L\infty (K_o;J) \colon
      \exists \, u_o \in \L\infty (\reali; J)
      \mbox{ with }
      \begin{array}{@{}r@{\,}c@{\,}l@{}}
        \rest{\ss_T u_o}{K_T}
        & =
        & u_T \,,
          \mbox{ and }
        \\
        \rest{u_o}{K_o}
        & =
        & \tilde u_o
      \end{array}
    \right\}.
  \end{equation}
\end{definition}

\noindent We regret that in the above notation
$\rest{\ii_T (u_T;J)}{K_o}$, the set $K_T$ is omitted for simplicity,
in spite of its relevance.

We now provide a simple specific extension of any map
$\widehat u_T \in \BV (K_T; J)$ to $u^*_T \in \BV (\reali; J)$ so that
the inverse design restricted to $K_o$ remains unaltered, provided
$K_o$ is the domain of dependency of $K_T$. More precisely:

\begin{theorem}
  \label{thm:dream1}
  Let $f$ satisfy~\ref{item:1}. Let $J$ be a non trivial closed real
  interval and $T$ be positive. Fix $x_1,x_2 \in \reali$ with
  $x_2 > x_1$ and choose a map $\widehat u_T \in \BV ([x_1,x_2];
  J)$. Define
  \begin{equation}
    \label{eq:22}
    u_T^* (x) = \left\{
      \begin{array}{l@{\qquad}l}
        \widehat u_T (x_1+)
        & x < x_1
        \\
        \widehat u_T (x)
        & x \in [x_1,x_2]
        \\
        \widehat u_T (x_2-)
        & x > x_2
      \end{array}
    \right.
  \end{equation}
  Then,
  \begin{displaymath}
    \rest{\ii_T (\widehat u_T;J)}{[\pi_{\widehat u_T} (x_1+), \pi_{\widehat u_T} (x_2-)]}
    =
    \rest{\ii_T (u^*_T;J)}{[\pi_{\widehat u_T} (x_1+), \pi_{\widehat u_T} (x_2-)]} \,.
  \end{displaymath}
\end{theorem}

In other words, setting in Theorem~\ref{thm:dream1} $K_T = [x_1,x_2]$
and
$K_o = \overline{\mathop{\rm co}} \; \pi_{\widehat u_T} (\pint{K}_T)$,
we have
\begin{displaymath}
  \rest{u^*_T}{K_T} = \widehat u_T
  \qquad \mbox{ and } \qquad
  \rest{\ii_T (\widehat u_T;J)}{K_o} =
  \rest{\ii_T(u^*_T;J)}{K_o} \,.
\end{displaymath}

In spite of its simplicity, the extension $u^*_T$ provided
in~\eqref{eq:22} of Theorem~\ref{thm:dream1} is sufficient to recover
the whole inverse design. The next result shows that any other
extension either gives the same result or gives the empty set.

\begin{theorem}
  \label{thm:dream2}
  Let $f$ satisfy~\ref{item:1}. Let $J$ be a non trivial closed real
  interval and $T$ be positive.  Fix $x_1,x_2 \in \reali$ with
  $x_2 > x_1$ and choose a map
  $\widehat u_T \in \L\infty ([x_1,x_2]; J)$ reachable at $t=T$ on
  $[x_1,x_2]$. Let $u_T \in \L\infty (\reali;J)$ be such that
  $\rest{u_T}{[x_1,x_2]} = \widehat u_T$. Then,
  \begin{description}
  \item[either:] $ \ii_T(u_T;J) = \emptyset$,
  \item[or:]
    $\rest{\ii_T (u_T;J)}{[\pi_{\widehat u_T} (x_1+), \pi_{\widehat
        u_T} (x_2-)]} = \rest{\ii_T (\widehat u_T;J)}{[\pi_{\widehat
        u_T} (x_1+), \pi_{\widehat u_T} (x_2-)]}$.
  \end{description}
\end{theorem}

Note that by Corollary~\ref{prop:Oleinik} the assumption that
$\widehat u_T$ be reachable ensures that $\widehat u_T$ has locally
bounded variation by~\ref{item:2}, hence its traces at $x_1$ and $x_2$
are well defined.

\begin{corollary}
  \label{cor:wakeup}
  Let $f$ satisfy~\ref{item:1}. Let $T$ be positive, $J$ be a non
  trivial closed real interval and fix $x_1$, $x_2$ in $\reali$ with
  $x_1 < x_2$.  Let $u_1,u_2 \in \L\infty (\reali; J)$ be such that
  $\rest{u_1}{[x_1,x_2]} = \rest{u_2}{[x_1,x_2]}$,
  $\ii_T (u_1) \neq \emptyset$ and $\ii_T (u_2) \neq \emptyset$. Then,
  with the notation~\eqref{eq:7} and~\eqref{eq:6}
  \begin{eqnarray}
    \label{eq:33}
    \rest{\ii_T (u_1;J)}{[\pi_{u_1} (x_1+), \pi_{u_1}
    (x_2-)]}
    & =
    & \rest{\ii_T (u_2;J)}{[\pi_{u_2} (x_1+), \pi_{u_2}
      (x_2-)]}\,,
    \\
    \label{eq:34}
    \rest{u^\flat_1}{[\pi_{u_1} (x_1+), \pi_{u_1}
    (x_2-)]}
    & =
    & \rest{u^\flat_2}{[\pi_{u_2} (x_1+), \pi_{u_2}
      (x_2-)]}\,,
    \\
    \label{eq:35}
    \rest{u^\sharp_1}{[\pi_{u_1} (x_1+), \pi_{u_1}
    (x_2-)]}
    & =
    & \rest{u^\sharp_2}{[\pi_{u_2} (x_1+), \pi_{u_2}
      (x_2-)]} \,.
  \end{eqnarray}
\end{corollary}

\noindent Equality~\eqref{eq:33} directly follows from
Theorem~\ref{thm:dream2}. Then, \eqref{eq:34} and~\eqref{eq:35} follow
combining~\eqref{eq:33} with Theorem~\ref{thm:1} setting for simplicity
$\widecheck x = x_1$.

\section{Application to Traffic}
\label{sec:traffic}

A macroscopic description of the flow of traffic along a highway tract
can be based on the well known Lighthill--Whitham and Richards
model~\cite{LighthillWhitham, Richards}, leading to the evolution
equation
\begin{equation}
  \label{eq:1}
  \partial_t \rho + \partial_x \left(\rho \, v (\rho)\right) = 0
  \qquad
  (t,x) \in \reali \times [0,L] \,,
\end{equation}
where $t$ is time, $x$ is the coordinate along the road,
$\rho = \rho (t,x)$ roughly measures the amount of vehicles per unit
length and $v = v (\rho)$ is the (mean) traffic speed corresponding to
the density $\rho$. As usual, we call $q = \rho \, v (\rho)$ the
vehicular flow. The space coordinate varies along the interval
$[0,L]$, with $L>0$. Note that~\eqref{eq:1} is neither a Cauchy
problem nor a standard initial - boundary value problem: nevertheless
it is a classical setup in traffic management.

At location $x=L$, the outflow $q_{out} = q_{out} (t)$ is
measured. The results in Section~\ref{subs:Notation} allow to exhibit
conditions implying that traffic underwent some critical event
(possibly an accident) and estimate where it happened. Moreover, they
also characterize vehicular traffic, providing properties and
constraints that any flow reconstruction in a data assimilation
procedure must enjoy or fulfill to be coherent with~\eqref{eq:1}. The
localization results in~\S~\ref{subsec:time-localization} allow us to
provide statements that are intrinsic to any (bounded) time interval
and that hold on the natural domain of dependency of the measured
data.

The speed law $v = v (\rho)$ plays in traffic modeling a role
analogous to that played in thermodynamics by the equation of
state. However, while equation of states can be rigorously justified
on the basis of physical assumptions, speed laws are typically
accepted or rejected on the basis of qualitative considerations. A
typical assumption on $v$ is
\begin{enumerate}[label=\bf(v)]
\item \label{item:6} $v \in \C2 ([0,R];\reali_+)$ is such that
  $v (R) = 0$ and
  $\frac{\d{}^2~}{\d\rho^2}\left(\rho \, v (\rho)\right) < 0$, for a
  fixed $R > 0$.
\end{enumerate}

A common problem in traffic modeling is the following: given the
traffic outflow measured at the position $x = L$, namely
$q_L (t) = \rho (t,L) \, v\left(\rho (t,L)\right)$, reconstruct the
function $\rho = \rho (t,x)$ for $x \in [0,L]$ assuming that the
outflow $q_L$ results from the -- unknown -- inflow at position $x=0$.

We are thus lead to exchange the roles of time $t$ and space $x$
in~\eqref{eq:1}, using as dependent variable the flow
$q (t,x) = \rho (t,x) \, v\left(\rho (t,x)\right)$ and refer to the
(backward) Cauchy problem
\begin{equation}
  \label{eq:2}
  \left\{
    \begin{array}{l}
      \partial_x q + \partial_t f (q) = 0
      \\
      q (t,L) = q_L (t)
    \end{array}
  \right.
  \qquad
  (t,x) \in \reali \times [0,L] \,.
\end{equation}
With reference to~\eqref{eq:1}, $f$ is (related to) the inverse of the
map $\rho \mapsto \rho \, v (\rho)$ on the interval where this map is
strictly increasing which, under assumption~\ref{item:6}, is the
interval $[0,\widehat q]$ where
$\widehat q = \max_{\rho \in [0,R]} \rho\,v (\rho)$. Note that the
choice of the congested interval where $\frac{\d{q}}{\d\rho} < 0$ is
not consistent with~\eqref{eq:2}.

\begin{proposition}
  \label{prop:equivalence}
  Let $v$ satisfy~\ref{item:6} and let $\widehat\rho$ be such that
  $\widehat\rho \, v (\widehat\rho) = \max_{[0,R]} \rho \, v (\rho)$.
  Define, for a fixed
  $\bar\rho \in \mathopen]0, \widehat\rho\mathclose[$ and for all
  $\rho \in [0, \bar\rho]$,
  \begin{displaymath}
    f (q) \colonequals \rho \iff q = \rho \, v (\rho) \,,
  \end{displaymath}
  Then:
  \begin{enumerate}[label=\bf(\arabic*)]
  \item \label{item:7} If, with reference to~\eqref{eq:1},
    $E \in \C2 ([0,\bar\rho]; \reali)$ is a convex entropy and $F$ a
    corresponding flux, then, $F\circ f$ is a convex entropy and
    $E \circ f$ is a corresponding flux for~\eqref{eq:2}.
  \item \label{item:8} If
    $\rho \in \L\infty (\reali \times [0,L]; [0,\bar\rho])$ is a weak
    solution to
    $\partial_t \rho + \partial_x \left(\rho \, v (\rho)\right) = 0$,
    then, the map
    $(t,x) \mapsto q (t,x) = \rho (t,x) \, v\left(\rho (t,x)\right)$
    is a weak solution to $\partial_x q + \partial_t f (q) = 0$
  \item \label{item:10} If in distributional sense
    \begin{equation}
      \label{eq:32}
      \partial_t E (\rho) + \partial_x F (\rho) \leq 0 \,,
    \end{equation}
    then, in distributional sense
    \begin{displaymath}
      \partial_x (F\circ f) (q) + \partial_t (E\circ f) (q) \leq 0 \,.
    \end{displaymath}
  \item \label{item:11} If~\eqref{eq:32} holds for any convex entropy,
    then, the trace $q_L$ defined by
    \begin{displaymath}
      \lim_{\delta \to 0+} \int_{-T}^T
      \modulo{q_L (t) - q (t,L-\delta)} \d{t} = 0
    \end{displaymath}
    for all $T > 0$, is well defined and $q$ is a weak solution
    to~\eqref{eq:2}.
  \end{enumerate}
\end{proposition}

\noindent The proofs of~\ref{item:7}, \ref{item:8} and~\ref{item:10}
are straightforward calculations, while~\ref{item:11} follows from the
regularity and convexity of $f$, thanks to~\cite{MR1869441}.

Condition~\ref{item:2} can then be interpreted as a minimal, necessary
but not sufficient, requirement for $q_L$ to be compatible with a
regular flow of traffic.

\begin{proposition}
  \label{prop:accidents}
  Under the assumptions of Proposition~\ref{prop:equivalence}, call
  $\bar q = \bar\rho \, v (\bar\rho)$. Fix a non trivial compact time
  interval $[T_1,T_2]$ and a measured traffic flow
  $q_{out} \in \L\infty ([T_1,T_2]; [0, \bar q])$. Then, $q_{out}$ is
  reachable at $x = L$, for some $L>0$, on $[T_1,T_2]$ in the sense of
  Definition~\ref{def:reachable} if and only if
  \begin{equation}
    \label{eq:4}
    \dfrac{1}{L}
    \geq
    \underset{\substack{T_1 \leq t_1 < t_2 \leq T_2}}{\esssup}
    \dfrac{f'\left(q_{out} (t_2)\right) - f'\left(q_{out} (t_1)\right)}{t_2-t_1} \,.
  \end{equation}
  Moreover, for such an $L$, define
  \begin{displaymath}
    \tau_1 = T_1 - L \; f'\left(q_{out} (T_1)\right)
    \quad \mbox{ and } \quad
    \tau_2 = T_2 - L \; f'\left(q_{out} (T_2)\right) \,.
  \end{displaymath}
  There exist
  $q^\flat, q^\sharp \in \L\infty([\tau_1,\tau_2];[0, \bar q]$ such
  that if $q_{in} \in \L\infty (\reali; [0, \bar q])$, the following
  statements are equivalent:
  \begin{enumerate}[label=\bf(\arabic*)]
  \item \label{item:12} There exists an entropy solution
    $q \in \L\infty (\reali \times [0,L];[0, \bar q])$ to
    \begin{equation}
      \label{eq:37}
      \partial_x q + \partial_t f (q) = 0
      \quad \mbox{ such that } \;
      \begin{array}{r@{\,}c@{\,}lr@{\,}c@{\,}l}
        q (t,0)
        & =
        & q_{in} (t)
        & \mbox{for } \; t
        & \in
        & [\tau_1,\tau_2] \,;
        \\
        q (t,L)
        & =
        & q_{out} (t)
        & \mbox{for } \; t
        & \in
        & [T_1,T_2] \,.
      \end{array}
    \end{equation}

  \item \label{item:9} For any $\tau \in [\tau_1,\tau_2]$
    \begin{equation}
      \label{eq:36}
      \int_{\tau_1}^{\tau}q_{in} (t) \d{t} \in
      \left[
        \int_{\tau_1}^{\tau}q^\flat (t) \d{t},
        \int_{\tau_1}^{\tau}q^\sharp (t) \d{t}
      \right] \,.
    \end{equation}
  \end{enumerate}
\end{proposition}

\noindent The characterization~\eqref{eq:4} is obtained applying
Corollary~\ref{prop:Oleinik}. Then, Theorem~\ref{thm:dream1},
Corollary~\ref{cor:wakeup} and Theorem~\ref{thm:1} allow to prove the
equivalence. Remark that Proposition~\ref{prop:accidents} is
\emph{intrinsic} to the time interval $[T_1,T_2]$, thanks in
particular to Corollary~\ref{cor:wakeup}.

\begin{remark}
  \label{rem:qflat}
  Localizing~\eqref{eq:7} in Theorem~\ref{thm:flat}, the function
  $q^\flat$, through one of its primitives $Q^\flat$, can be computed
  from the measured data $q_{out}$ for any $\tau \in [\tau_1,\tau_2]$
  \begin{displaymath}
    Q^\flat (\tau)
    \colonequals
    \sup_{t \in [T_1,T_2]} \left[
      \int_{T_1}^t q_{out} (s) \d{s}
      -
      L \, f^*\left(\frac{t-\tau}{L}\right)
    \right] \,.
  \end{displaymath}
\end{remark}

From the traffic management point of view, as soon as
condition~\eqref{eq:4} is violated, one can infer that the standard
flow of traffic was altered -- possibly by an accident -- at a
distance $L$ from the measuring site. For such $L$, a time $\tau$ at
which traffic resumes after the road was blocked, is a time where
$Q^\flat$ is not differentiable.

\section{Technical Details}
\label{sec:Tech}

Recall first the following elementary definitions and properties.

\begin{lemma}
  \label{lem:SC0} Let $f \colon \reali \to \reali$ be strongly
  convex. Then, for all $A \in \reali$, there exists
  $\alpha \in \reali$ such that for all $x \in \reali$,
  $f (x) \geq \alpha + A\,\modulo{x} $.
\end{lemma}

\begin{definition}
  \label{def:Legendre}
  Let $f \colon \reali \to \reali$ be convex. Its \emph{Legendre
    Transform} is the map $f^* \colon \reali \to \reali$ defined by
  $f^* (y) \colonequals \sup_{x \in \reali} \left(y \, x - f
    (x)\right)$. for all $y \in \reali$.
\end{definition}

Note that for any compact interval $K$, the values of $f^*$ on
$f' (K)$ depend exclusively on the restriction of $f$ to $K$.

\begin{lemma}[{\cite[Theorem~3, \S~3.3.2]{Evans}}]
  \label{lem:SC}
  If $f \colon \reali \to \reali$ is strongly convex and $f^*$ is its
  Legendre transform, then
  \begin{enumerate}[label={(L\arabic*)}]
  \item \label{item:3} for all $y \in \reali$,
    $f^* (y) = y \, (f')^{-1} (y) - f\left((f')^{-1} (y)\right)$;
  \item \label{item:4} for all $y \in \reali$,
    $(f^*)' (y) = (f')^{-1} (y)$;
  \item \label{item:5} $f^*$ is strongly convex in the sense it
    satisfies~\ref{item:1}.
  \end{enumerate}
\end{lemma}

\subsection{Proofs Related to \S~\ref{subsec:inverse-design}}
\label{subsec:proofs-inverse-design}

\begin{proofof}{Theorem~\ref{thm:flat}}
  The proof is divided into several short steps. Using the
  notation~\eqref{eq:7}, define the set valued map
  \begin{equation}
    \label{eq:8}
    \begin{array}{ccccc}
      \mathcal{M}
      & \colon
      & \reali
      & \to
      & \mathcal{P} (\reali)
      \\
      &
      & x
      & \mapsto
      &\left\{
        \xi \in \reali \colon
        U_T (\xi) = U_o^\flat (x) + T \, f^*\!\left(\frac{\xi - x}{T}\right)
        \right\}
    \end{array}
  \end{equation}

  \paragraph{1. $U_o^\flat \in \Lip (\reali; \reali)$ and for
    a.e.~$x \in \reali$, $\frac{\d{~}}{\d x}U_o^\flat (x) \in J$.}

  A direct consequence of the assumptions on $u_T$ is that, for $h$
  positive, by~\eqref{eq:7}
  \begin{equation}
    \label{eq:9}
    U_T (x) - h \, \essinf_{\reali} u_T
    \geq
    U_T (x-h)
    \geq
    U_T (x) - h \, \esssup_{\reali} u_T \,.
  \end{equation}
  Moreover, for any $x_1, x_2 \in \reali$ with $x_1 < x_2$, setting
  $\xi = x - (x_2 - x_1)$,
  \begin{displaymath}
    U_o^\flat (x_1)
    =
    \sup_{\xi \in \reali}
    \left[
      U_T (\xi) - T \, f^*\left(\frac{\xi-x}{T}\right)
    \right]
    =
    \sup_{x \in \reali}
    \left[
      U_T \left(x - (x_2-x_1)\right)
      - T \, f^*\left(\frac{x-x_2}{T}\right)
    \right]
  \end{displaymath}
  so that by~\eqref{eq:9} with $h = x_2-x_1$, we have the two
  estimates
  \begin{eqnarray*}
    U_o^\flat (x_1)
    & \leq
    & \sup_{x \in \reali}
      \left(U_T (x) - T \, f^*\left(\frac{x-x_2}{T}\right)\right)
      - (x_2 - x_1) \, \essinf_{\reali} u_T
    \\
    & =
    & U_o^\flat (x_2)
      - (x_2 - x_1) \, \essinf_{\reali} u_T \,;
    \\
    U_o^\flat (x_1)
    & \geq
    & \sup_{x \in \reali}
      \left(U_T (x) - T \, f^*\left(\frac{x-x_2}{T}\right)\right)
      - (x_2-x_1) \, \esssup_{\reali} u_T
    \\
    & =
    & U_o^\flat (x_2) - (x_2-x_1) \, \esssup_{\reali} u_T \,,
  \end{eqnarray*}
  proving the Lipschitz continuity of $U_o^\flat$ and, since
  $[\essinf_{\reali} u_T, \esssup_{\reali} u_T] \subseteq J$,
  completing the proof of the claim, thanks to Rademacher's
  Theorem~\cite[Theorem~6, \S~5.8.3]{Evans}.\claimend

  \paragraph{2. For all $x$, the set $\mathcal{M} (x)$ is not empty.}
  (In other words, the $\sup$ in~\eqref{eq:7} is a maximum).
  By~\ref{item:1} and Lemma~\ref{lem:SC}, the map $f^*$ is strongly
  convex, the map $U_T$ is globally Lipschitz continuous, and so
  sublinear at $\pm\infty$, hence by Lemma~\ref{lem:SC0}
  $\lim_{\xi \to -\infty} \left(U_T (\xi) - T \,
    f^*\!\left(\frac{\xi-x}{T}\right)\right) = -\infty$ and
  $\lim_{\xi \to +\infty} \left(U_T (\xi) - T \,
    f^*\!\left(\frac{\xi-x}{T}\right)\right) = -\infty$ for all
  $x \in \reali$. An application of {Weierstrass} Theorem completes
  the proof of the claim.\claimend

  \paragraph{3. There exists $R>0$ such that if $x \in \reali$ and
    $\xi \in \mathcal{M} (x)$, then $\modulo{x-\xi} \leq R$.} Let
  $\kappa$ be a Lipschitz constant for $U_T$. By~\ref{item:1} and
  Lemma~\ref{lem:SC}, $f^*$ is strongly convex, so that we can apply
  Lemma~\ref{lem:SC0} to $f^*$ with $A = \kappa + 1$. Hence,
  \begin{eqnarray*}
    U_T (\xi) - T \, f^*\left(\frac{\xi-x}{T}\right)
    & \leq
    & U_T (x)
      +
      \kappa \, \modulo{\xi-x}
      -
      T \, \alpha - (\kappa+1) \, \modulo{\xi-x}
    \\
    & \leq
    & U_T (x) - T f^* (0) + T \left(f^* (0) - \alpha\right) - \modulo{\xi-x}
  \end{eqnarray*}
  and the choice of $\xi$ as a maximizer for the left hand side above
  ensures that
  $\modulo{x-\xi} \leq T \left(f^* (0) - \alpha\right)$.\claimend

  \paragraph{4. $\mathcal{M}$ is monotone increasing.} By this, we
  mean that if $x_1, x_2 \in \reali$ and $x_1 < x_2$, then for all
  $\xi_1 \in \mathcal{M} (x_1)$ and for all
  $\xi_2 \in \mathcal{M} (x_2)$, it holds that $\xi_1 \leq \xi_2$.

  Proceed by contradiction and assume that $x_1 < x_2$ but
  $\xi_1 > \xi_2$. Define
  \begin{displaymath}
    A = \dfrac{\xi_1 - x_1}{T} \,,\quad
    B = \dfrac{\xi_2 - x_2}{T} \,,\quad
    C = \dfrac{\xi_2 - x_1}{T} \; \mbox{ and } \;
    D = \dfrac{\xi_1 - x_2}{T} \,.
  \end{displaymath}
  Clearly, $A+B = C+D$ and $A>C>B$, $A>D>B$.  By construction, there
  exists a (unique) $\theta \in \left]0, 1 \right[$ such that
  \begin{displaymath}
    C = \theta \, A + (1-\theta) B
    \quad \mbox{ and } \quad
    D = (1-\theta) A + \theta \, B \,.
  \end{displaymath}
  The strong convexity of $f^*$ then ensures that
  \begin{displaymath}
    f^* (C) + f^* (D) < f^* (A) + f^* (B) \,.
  \end{displaymath}
  The above choices of $x_1,x_2,\xi_1,\xi_2$ imply that
  \begin{displaymath}
    U_o^\flat (x_1) = U_T (\xi_1) - T \, f^* (A)
    \quad \mbox{ and } \quad
    U_o^\flat (x_2) = U_T (\xi_2) - T \, f^*(B) \,,
  \end{displaymath}
  so that
  \begin{eqnarray*}
    U_o^\flat (x_1) + U_o^\flat (x_2)
    & <
    & U_T (\xi_1) + U_T (\xi_2) - T \, f^* (C) - T \, f^* (D)
  \end{eqnarray*}
  which in turn implies that at least one of the following
  inequalities hold:
  \begin{displaymath}
    U_o^\flat (x_1) < U_T (\xi_2) - T \, f^* (C)
    \quad \mbox{ or } \quad
    U_o^\flat (x_2) < U_T (\xi_1) - T \, f^* (D) \,.
  \end{displaymath}
  Both inequalities above contradict the definition~\eqref{eq:7} of
  $U_o^\flat$.\claimend

  \paragraph{5. For all $x$, the set $\mathcal{M} (x)$ is a compact
    interval.}

  By Claim~3.~above, for all $x \in \reali$, the set $\mathcal{M} (x)$
  is bounded since $\mathcal{M} (x) \subseteq [x-R, x+R]$. The
  definition~\eqref{eq:8} of $\mathcal{M}$ shows that for every fixed
  $x$, the set $\mathcal{M} (x)$ is closed, thanks to the regularity
  of $U_o^\flat$ proved in Claim~1 and that of $U_T$ and $f^*$.

  Fix $x \in \reali$. To prove that $\mathcal{M} (x)$ is an interval,
  choose $\xi_1, \xi_2 \in \mathcal{M} (x)$ with $\xi_1 < \xi_2$ and
  introduce the Lipschitz continuous map
  $\phi \colon [\xi_1, \xi_2] \to \reali$ by
  $\phi (\xi) \colonequals U_T (\xi)- T
  f^*\left((\xi-x)/T\right)$. Proceed by contradiction and.
  by~\eqref{eq:7} and~\eqref{eq:8}, assume that there exists a
  $\xi \in \mathopen]\xi_1, \xi_2 \mathclose[$ such that
  $\phi (\xi) < \phi (\xi_1) = \phi (\xi_2) = \max \phi$. Then, the
  Lipschitz continuity of $\phi$ ensures that $\phi$ is differentiable
  a.e.~and
  \begin{displaymath}
    \begin{array}{@{\!}rl@{\quad}|@{\quad}rl@{\!}}
      &0 > \phi (\xi) - \phi (\xi_1) = \int_{\xi_1}^\xi \phi' (y) \d{y}
      &
      &0 > \phi (\xi_2) - \phi (\xi) = \int_\xi^{\xi_2} \phi' (y) \d{y}
      \\
      \Rightarrow
      &\exists \, y_1 \in \left]\xi_1, \xi\right[ \colon
        \left\{
        \begin{array}{@{}l@{}}
          \phi \mbox{ is differentiable at } y_1
          \\
          \mbox{and }  \phi' (y_1) < 0
        \end{array}
      \right.
      &  \Rightarrow
      &\exists \, y_2 \in \left]\xi, \xi_2\right[ \colon
        \left\{
        \begin{array}{@{}l@{}}
          \phi \mbox{ is differentiable at } y_2
          \\
          \mbox{and }  \phi' (y_2) > 0
        \end{array}
      \right.
      \\
      \Rightarrow
      & U_T' (y_1) < (f^*)'\left((y_1-x)/T\right)
      & \Rightarrow
      & U_T' (y_2) > (f^*)'\left((y_2-x)/T\right)
      \\
      \Rightarrow
      & f'\left(U_T' (y_1)\right) < (y_1-x)/T
      & \Rightarrow
      & f'\left(U_T' (y_2)\right) > (y_2-x)/T
    \end{array}
  \end{displaymath}
  so that
  \begin{displaymath}
    f'\left(U_T' (y_2)\right) - f'\left(U_T' (y_1)\right)
    >
    \dfrac{y_2-y_1}{T}
  \end{displaymath}
  which contradicts that $u_T$ satisfies Condition~\ref{item:2} at
  $T$ since $y_1 < \xi < y_2$.\claimend

  \paragraph{6. $\mathcal{M}$ is surjective, in the sense that
    $\bigcup_{x \in \reali} \mathcal{M} (x) = \reali$.}
  Fix $\xi$ in $\reali$. We seek an $x \in \reali$ such that
  $\xi = \mathcal{M} (x)$. Thanks to Claim~5., it is sufficient to
  prove that
  $\min \mathcal{M} (x) \leq \xi \leq \max \mathcal{M} (x)$.

  To this aim, define
  $x_* = \sup \left\{x \in \reali \colon \min \mathcal{M} (x) <
    \xi\right\}$ and note that
  $\min \mathcal{M} (x_*- \frac{1}{n}) < \xi$. The sequence $\xi_n$
  defined by $\xi_n = \min \mathcal{M} (x_* - \frac1n)$ is increasing
  by Claim~4.~above and it is bounded by the above choices. Thus,
  there exits a real $\xi_*$ such that
  $\lim_{n\to +\infty} \xi_n = \xi_*$. Definition~\eqref{eq:8} and
  Claim~1.~then ensure that $\xi_* \in \mathcal{M} (x_*)$ and
  $\xi_* \leq \xi$. As a consequence,
  $\min \mathcal{M} (x_*) \leq \xi$.

  To prove the other bound $\max \mathcal{M} (x_*) \geq \xi$, proceed
  by contradiction. Assume that $\max \mathcal{M}(x_*) < \xi$ and
  consider the sequence $y_n = \max \mathcal{M} (x_*+\frac{1}{n})$,
  which is decreasing and bounded below, so that it admits a limit
  $y_*$, with $y_* \in \mathcal{M} (x_*)$. Thus, $y_* < \xi$ and for
  at least one index $n$, we have $y_n < \xi$, so that
  $\min \mathcal{M} (x_*+\frac{1}{n}) \leq y_n < \xi$ which
  contradicts the inequality $x_*+\frac{1}{n} > x_*$.\claimend

  \paragraph{7. Conclusion.} By~\eqref{eq:7}, for all
  $x,\xi \in \reali$, we have that
  $U_T (\xi) \leq U_o^\flat (x) + T \,
  f^*\left(\frac{\xi-x}{T}\right)$. By Claim~6., for all
  $\xi \in \reali$, there exists $x \in \reali$ such that
  $U_T (\xi) = U_o^\flat (x) + T \,
  f^*\left(\frac{\xi-x}{T}\right)$. Therefore,
  \begin{displaymath}
    U_T (\xi)
    =
    \inf_{x \in \reali} U_o^\flat (x)
    +
    T \, f^*\left(\dfrac{\xi-x}{T}\right)
  \end{displaymath}
  proving, by Hopf-Lax Formula\cite[Theorem~4, Section~3,
  Chapter~3]{Evans}, that $U_o^\flat \in \II_T (U_T; J)$ and, hence,
  that $u_o^\flat \in \ii_T (u_T; J)$ by~\eqref{eq:10}.
\end{proofof}

\begin{proofof}{Proposition~\ref{prop:reverse}}
  Call $g (u) \colonequals f (-u)$. By Definition~\ref{def:Legendre},
  $g^* (p) = f^* (-p)$ for all $p \in \reali$.
  Definition~\eqref{eq:7} of $U_o^\flat$ implies that for all
  $x \in \reali$
  \begin{displaymath}
    -U_o^\flat (x)
    =
    \inf_{\xi \in \reali} \left(
      -U_T (\xi) + T \; f^*\left(-\frac{x-\xi}{T}\right)
    \right)
    =
    \inf_{\xi \in \reali} \left(
      -U_T (\xi) + T \; g^*\left(\frac{x-\xi}{T}\right)
    \right)
  \end{displaymath}
  By the Hopf-Lax Formula\cite[Theorem~4, Section~3,
  Chapter~3]{Evans}, we have $V (T,x) = - U_o^\flat (x)$. An
  application of~\cite[Theorem~1.1]{MR1923522} as in~\eqref{eq:10}
  ensures that $v (T,x) = - u_o^\flat (x)$ and that
  $v (t,x) = \partial_x V (t,x)$. Finally, the desired bound on the
  total variation follows
  from~\cite[Theorem~6.1]{BressanLectureNotes}.
\end{proofof}

\begin{proofof}{Corollary~\ref{cor:clco}}
  The equality of the total variations follows from
  Proposition~\ref{prop:reverse}
  and~\cite[Theorem~6.1]{BressanLectureNotes}. Then, left continuous
  representatives can be exploited. Finally, the equality
  $\ss_Tu_o^\flat = u_T$ implies that
  $\overline{\mathop{\rm co}} \; u_T (\reali) \subseteq
  \overline{\mathop{\rm co}} \; u_o^\flat (\reali)$. The other
  inclusion follows from Theorem~\ref{thm:flat} choosing
  $J \colonequals \overline{\mathop{\rm co}} \; u_T (\reali)$.
\end{proofof}

\begin{proofof}{Proposition~\ref{prop:characterization}}
  The leftmost equivalence in~\eqref{eq:3} follows directly
  from~\eqref{eq:10}.

  Assume now that $U_o \in \II_T (U_T; J)$. Then, to prove~$(i)$, start
  from~\eqref{eq:17}:
  \begin{displaymath}
    \begin{array}{cl@{\qquad}r@{\,}c@{\,}l}
      & \forall\, x \in \reali
      & U_T (x)
      & =
      & \inf_{\xi \in \reali} U_o (\xi) + T \, f^*\left(\frac{x-\xi}{T}\right)
      \\
      \implies
      & \forall\, x,\xi \in \reali
      & U_T (x)
      & \leq
      & U_o (\xi) + T \, f^*\left(\frac{x-\xi}{T}\right)
      \\
      \implies
      & \forall\, x,\xi \in \reali
      & U_o (\xi)
      & \geq
      & U_T (x) - T \, f^*\left(\frac{x-\xi}{T}\right)
      \\
      \implies
      & \forall\, \xi \in \reali
      & U_o (\xi)
      & \geq
      & \sup_{x \in \reali} U_T (x) - T \, f^*\left(\frac{x-\xi}{T}\right)
        = U_o^\flat (\xi)
    \end{array}
  \end{displaymath}
  by~\eqref{eq:7}. Condition~$(ii)$ follows from
  Proposition~\ref{thm:two_min_renewed} and~$(iii)$ is a consequence
  of~\eqref{eq:5}.

  Finally, assume that~$(i)$, $(ii)$ and~$(iii)$ hold. Define
  \begin{equation}
    \label{eq:20}
    \forall\,x \in \reali \qquad
    \widehat U_T (x)
    \colonequals
    \inf_{\xi \in \reali} U_o (\xi) + T \, f^*\left(\frac{x-\xi}{T}\right)
  \end{equation}
  so that $\widehat U_T\geq U_T$ by~$(i)$ and since
  $U_o^\flat \in \II_T(U_T;J)$. On the other hand, for any
  $x \in \reali$,
  \begin{flalign*}
    U_T (x) \; = \; %
    & U_o^\flat\left(\pi_{u_T} (x)\right) + T \;
    f^*\left(\frac{x-\pi_{u_T} (x)}{T}\right) %
    & [\mbox{By Proposition~\ref{thm:two_min_renewed}}]
    \\
    \; = \; %
    & U_o\left(\pi_{u_T} (x)\right) + T \; f^*\left(\frac{x-\pi_{u_T}
        (x)}{T}\right) %
    & [\mbox{By~$(ii)$]}
    \\
    \; \geq \; %
    & \widehat U_T (x) %
    & [\mbox{By~\eqref{eq:20}]}
  \end{flalign*}
  completing the proof, thanks to~\eqref{eq:17}, $(iii)$
  and~\eqref{eq:5}.
\end{proofof}

\begin{proofof}{Proposition~\ref{prop:sharp}}
  Note first that $U_o^\sharp$ is well defined. Indeed, fix any
  $\bar x \in \overline{\pi_{U_T'} (\reali)}$. Then, for all
  $U_o \in \II_T (U_T; J)$, we have
  $U_o (\bar x) = U_o^\flat (\bar x)$, by~$(i)$ in
  Proposition~\ref{prop:characterization}. Then, calling
  $M \colonequals \max \left\{\modulo{w} \colon w \in J\right\}$, for
  all $x \in \reali$, we have that for all $U_o \in \II_T (U_T; J)$,
  $\modulo{U_o (x)} \leq \modulo{U_o^\flat (\bar x)} + M
  \modulo{x-\bar x}$, showing that the $\sup$ in~\eqref{eq:6} is in
  $\reali$.

  We now prove that $U_o^\sharp$ satisfies~$(iii)$. Indeed, fix
  $x_1,x_2 \in \reali$ with $x_1 < x_2$. Then, for any
  $U_o \in \II_T (U_T; J)$
  \begin{equation}
    \label{eq:12}
    \min J \leq \dfrac{U_o (x_2) - U_o (x_1)}{x_2 - x_1} \leq \max J \,.
  \end{equation}
  So that
  \begin{displaymath}
    \begin{array}{rcl@{\qquad}|@{\qquad}rcl}
      U_o (x_2)
      & \leq
      & U_o (x_1) + (x_2-x_1) \max J
      & U_o (x_1)
      & \leq
      & U_o (x_2) - (x_2-x_1) \min J
      \\
      U_o (x_2)
      & \leq
      & U_o^\sharp (x_1) + (x_2-x_1) \max J
      & U_o (x_1)
      & \leq
      & U_o^\sharp (x_2) - (x_2-x_1) \min J
      \\
      U_o^\sharp (x_2)
      & \leq
      & U_o^\sharp (x_1) + (x_2-x_1) \max J
      & U_o^\sharp (x_1)
      & \leq
      & U_o^\sharp (x_2) - (x_2-x_1) \min J
    \end{array}
  \end{displaymath}
  proving~$(iii)$.

  To complete the proof, simply observe that $U_o^\sharp$ satisfies
  ~$(i)$ and $(ii)$ by construction.
\end{proofof}

\begin{proofof}{Theorem~\ref{thm:1}}
  To prove the equality in~\eqref{eq:14}, note that the inclusion
  $\subseteq$ follows from~$(i)$--$(iii)$ in~\eqref{eq:3} and from the
  definition of $U_o^\sharp$ in~\eqref{eq:6}. On the other hand, $(i)$
  and~$(iii)$ are obvious, while Proposition~\ref{prop:sharp}
  ensures~$(ii)$ in~\eqref{eq:3}, which then implies the other
  inclusion~$\supseteq$. Convexity is now straightforward. Compactness
  in the said topology follows from~\eqref{eq:14} and Ascoli--Arzel\`a
  Theorem~\cite[\S~C.7]{Evans}, which can be applied thanks to the
  compactness of $J$.

  The correspondences described in~\eqref{eq:10} now ensure the
  equality
  \begin{displaymath}
    \ii_T (u_T;J)
    =
    \left\{
      u_o \in \L\infty (\reali; J)
      \colon
      \exists \, U_o \in \II_T (U_T;J)
      \mbox{ such that }
      u_o = U_o'
    \right\}
  \end{displaymath}
  and therefore the convexity of $\ii_T (u_T;J)$. To
  prove~\eqref{eq:16}, use~\eqref{eq:14} and~\eqref{eq:11} and recall
  that we set $\widecheck y = \pi_{u_T} (\widecheck x)$. Sequential
  compactness follows from~\eqref{eq:16} and from the boundedness of
  $\ii_T (u_T;J)$ in the weak-$*$ topology, see~\cite[(ii) in
  \S~4.3.C]{MR2759829}. 
\end{proofof}

\subsection{Proofs Related to \S~\ref{subsec:time-localization}}
\label{subsec:proofs-time-localization}

\begin{lemma}
  \label{lem:glue}
  Let $f$ satisfy~\ref{item:1}. Fix
  $u_o, v_o \in \L\infty (\reali; \reali)$ and call $u$, respectively,
  $v$ the weak entropy solution to
  \begin{displaymath}
    \left\{
      \begin{array}{l}
        \partial_t u + \partial_x f (u) = 0
        \\
        u (0,x) = u_o (x)\,,
      \end{array}
    \right.
    \quad \mbox{ respectively } \quad
    \left\{
      \begin{array}{l}
        \partial_t v + \partial_x f (v) = 0
        \\
        v (0,x) = v_o (x) \,.
      \end{array}
    \right.
  \end{displaymath}
  Assume there exist $T>0$ and $\bar x \in \reali$ such that
  \begin{equation}
    \label{eq:24}
    u (T, \bar x-) = v (T, \bar x-) \,.
  \end{equation}
  Then, the map
  \begin{displaymath}
    w (t,x)
    \colonequals
    \left\{
      \begin{array}{l@{\qquad}r@{\,}c@{}l}
        u (t,x)
        & x
        & <
        & \pi (t)
        \\
        v (t,x)
        & x
        & \geq
        & \pi (t)
      \end{array}
    \right.
    \quad \mbox{ where } \quad
    \pi (t) \colonequals \bar x + (t-T) \; f'\left(u (T, \bar x-)\right)
  \end{displaymath}
  is a weak entropy solution to
  \begin{displaymath}
    \left\{
      \begin{array}{l}
        \partial_t w + \partial_x f (w) = 0
        \\
        w (0,x) = w_o (x)
      \end{array}
    \right.
    \quad \mbox{ where } \quad
    w_o (x)
    \colonequals
    \left\{
      \begin{array}{l@{\qquad}r@{\,}c@{}l}
        u_o (x)
        & x
        & <
        & \pi (0)
        \\
        v_o (x)
        & x
        & \geq
        & \pi (0) \,.
      \end{array}
    \right.
  \end{displaymath}
  An entirely similar statement holds replacing the left trace with
  the right trace in~\eqref{eq:24} and in the definition of $\pi$.
\end{lemma}

Remark that, since $T > 0$, the existence of the traces
in~\eqref{eq:24} is ensured by the uniform convexity of $f$
by~\cite[Theorem~11.2.2]{MR3468916}. Moreover, by~\eqref{eq:24}, the
line $x = \pi (t)$ is a minimal backward
characteristics~\cite[Theorem~11.1.3]{MR3468916} common to both $u$
and $v$.

\begin{proofof}{Lemma~\ref{lem:glue}}
  Define
  \begin{displaymath}
    \begin{array}{@{}r@{\,}c@{\,}l@{;\quad}r@{\,}c@{\,}l@{;\quad}r@{\,}c@{\,}l@{}}
      U_o (x)
      & \colonequals
      & \int_{\pi (0)}^x u_o (\xi) \d\xi
      & U (t,x)
      & \colonequals
      & \int_{\pi (t)}^x u (t,\xi) \d\xi
      & p_u (t,x)
      & \colonequals
      & x - t \, f'\left(u (t, x-)\right);
      \\
      V_o (x)
      & \colonequals
      & \int_{\pi (0)}^x v_o (\xi) \d\xi
      & V (t,x)
      & \colonequals
      & \int_{\pi (t)}^x v (t,\xi) \d\xi
      & p_v (t,x)
      & \colonequals
      & x - t \, f'\left(v (t, x-)\right);
      \\
      W_o (x)
      & \colonequals
      & \int_{\pi (0)}^x w_o (\xi) \d\xi
      & W (t,x)
      & \colonequals
      & \int_{\pi (t)}^x w (t,\xi) \d\xi
      & p_w (t,x)
      & \colonequals
      & x - t \, f'\left(w (t, x-)\right).
    \end{array}
  \end{displaymath}
  By~\eqref{eq:24}, $\pi$ is a genuine generalized characteristic of
  both $u$ and $v$, and since genuine backward characteristics do not
  cross~\cite[Corollary~11.1.2]{MR3468916}, we have
  \begin{displaymath}
    \forall\, (t,x) \in [0,T] \times \reali
    \qquad
    \begin{array}{r@{\,}c@{\,}l@{\implies}r@{\,}c@{\,}l}
      x
      & <
      & \pi (t)
      & p_u (t,x)
      & =
      & p_w (t,x)
      \\
      x
      & \geq
      & \pi (t)
      & p_v (t,x)
      & =
      & p_w (t,x)
    \end{array}
  \end{displaymath}
  so that, by~\cite[Theorem~11.4.3]{MR3468916},
  \begin{displaymath}
    \begin{array}{r@{\,}c@{\,}l@{\quad}|@{\quad}r@{\,}c@{\,}l}
      &
      & \mbox{if }x < \pi (t)\,{:}
      &
      &
      & \mbox{if }x > \pi (t)\,{:}
      \\
      W (t,x)
      & =
      & U (t,x)
      & W (t,x)
      & =
      & V (t,x)
      \\
      & =
      & U_o\left(p_u (t,x)\right) + t \, f^*\left(\frac{x-p_u (t,x)}{t}\right)
      &
      & =
      & V_o\left(p_v (t,x)\right) + t \, f^*\left(\frac{x-p_v (t,x)}{t}\right)
      \\
      & =
      & W_o\left(p_w (t,x)\right) + t \, f^*\left(\frac{x-p_w (t,x)}{t}\right)
        \,;
      &
      & =
      & W_o\left(p_w (t,x)\right) + t \, f^*\left(\frac{x-p_w (t,x)}{t}\right)
        \,.
    \end{array}
  \end{displaymath}
  Hence, by~\cite[Theorem~4, Section~3, Chapter~3]{Evans}
  and~\cite[Theorem~11.4.3]{MR3468916}, $W$ is the viscosity solution
  to
  \begin{displaymath}
    \left\{
      \begin{array}{l}
        \partial_t W + f (\partial_x W) = 0
        \\
        W (0,x) = W_o (x)
      \end{array}
    \right.
  \end{displaymath}
  which yields the proof by the correspondence
  in~\cite[Theorem~1.1]{MR1923522}, see also~\eqref{eq:10}.
\end{proofof}

\begin{lemma}
  \label{lem:TBA}
  Let $f$ satisfy~\ref{item:1}. Let $J$ be a non trivial closed real
  interval and $T$ be positive. Fix $x_1,x_2 \in \reali$ with
  $x_2 > x_1$ and $u_o,u_T \in \L\infty (\reali;J)$ such that
  \begin{equation}
    \label{eq:23}
    u_T = \ss_T u_o \,,
  \end{equation}
  Define, as shown in Figure~\ref{fig:TBA},
  \begin{eqnarray*}
    u_o^* (x)
    & \colonequals
    & \left\{
      \begin{array}{@{\,}l@{\qquad}l}
        u_o \left(\pi_{u_T}(x_1+)+\right)
        & x < \pi_{u_T} (x_1+)
        \\
        u_o (x)
        & x \in [\pi_{u_T} (x_1+),\pi_{u_T}(x_2-)]
        \\
        u_o \left(\pi_{u_T}(x_2-)-\right)
        & x > \pi_{u_T}(x_2-)
      \end{array}
          \right.
    \\
    u_T^* (x)
    & \colonequals
    & \left\{
      \begin{array}{@{\,}l@{\qquad}l}
        u_T (x_1+)
        & x < x_1
        \\
        u_T (x)
        & x \in [x_1, x_2]
        \\
        u_T (x_2-)
        & x > x_2 \,.
      \end{array}
          \right.
  \end{eqnarray*}
  Then,
  \begin{displaymath}
    u_T^* = \ss_T u_o^* \,.
  \end{displaymath}
\end{lemma}

\begin{figure}[!h]
  \begin{center}
    \begin{tikzpicture}[baseline]
      \draw[thick,-stealth] (-5,0) -- (10,0); %
      \draw[thick,-stealth] (0,-1) -- (0,5); %
      \draw (-0.5,5) node {$t$}; %
      \draw (10,-0.5) node {$x$}; %
      \draw (-0.5,3.75) node {$T$}; %
      \draw (-0.5, -0.5) node {$0$}; %
      \draw[dashed] (-5,4) -- (10,4); %
      \draw (-2.5,0) -- (-1,4); %
      \draw (-3, -0.5) node {$\pi_{u_T} (x_1+)$}; %
      \draw (-1,4.5) node {$x_1$}; %
      \draw (7,0) -- (5.5,4); %
      \draw (6.5, -0.5) node {$\pi_{u_T} (x_2-)$}; %
      \draw (5.5,4.5) node {$x_2$}; %
      \draw (2.5,0.25) node {$u_o^* (x) {=} u_o (x)$}; %
      \draw (2.25,4.25) node {$\ss_T u_o^* (x) {=} u_T (x)$}; %
      \draw (8.25, 2) node {$\ss_t u_o^* (x) {=} u_T (x_2-)$}; %
      \draw (-3.5, 2) node {$\ss_t u_o^* (x) {=} u_T (x_1+)$};
    \end{tikzpicture}
  \end{center}
  \caption{Notations used in Lemma~\ref{lem:TBA}.}
  \label{fig:TBA}
\end{figure}

\begin{proofof}{Lemma~\ref{lem:TBA}}
  Clearly, both constant maps $(t,x) \mapsto u_T (t_1+)$ and
  $(t,x) \mapsto u_T (t_2-)$ are entropy solutions to
  $\partial_t u + \partial_x f (u) = 0$. Thus, the proof follows by an
  application of Lemma~\ref{lem:glue}.
\end{proofof}

\begin{proofof}{Theorem~\ref{thm:dream1}}
  Definition~\eqref{eq:30}, implies the inclusion
  \begin{displaymath}
    \underbrace{\rest{\ii_T (\widehat u_T;J)}{[\pi_{\widehat u_T} (x_1+),
        \pi_{\widehat u_T} (x_2-)]}}_{\mbox{here } K_T = [x_1,x_2]}
    \supseteq
    \underbrace{\rest{\ii_T (u^*_T;J)}{[\pi_{\widehat u_T} (x_1+), \pi_{\widehat u_T} (x_2-)]}}_{\mbox{here } K_T = \reali} \,.
  \end{displaymath}
  If the set in the left hand side above, namely
  $\rest{\ii_T (\widehat u_T;J)}{[\pi_{\widehat u_T} (x_1+),
    \pi_{\widehat u_T} (x_2-)]}$ is empty, then the proof trivially
  follows. Otherwise, if
  $\rest{\ii_T (\widehat u_T;J)}{[\pi_{\widehat u_T} (x_1+),
    \pi_{\widehat u_T} (x_2-)]} \neq \emptyset$, choose any
  $\widetilde u_o \in \rest{\ii_T (\widehat u_T;J)}{[\pi_{\widehat
      u_T} (x_1+), \pi_{\widehat u_T} (x_2-)]}$. Then, with reference
  to~\eqref{eq:30}, there exists a $u_o \in \L\infty (\reali; J)$ such
  that
  $\rest{u_o}{[\pi_{\widehat u_T} (x_1+), \pi_{\widehat u_T} (x_2-)]}
  = \widetilde u_o$ and
  $\rest{\ss_T {u_o}}{[x_1, x_2]} = \widehat u_T$. An application of
  Lemma~\ref{lem:TBA} completes the proof.
\end{proofof}

\begin{proofof}{Theorem~\ref{thm:dream2}}
  Assume that $\ii_T (u_T;J) \neq
  \emptyset$. Definition~\eqref{eq:30}, ensures the inclusion
  \begin{displaymath}
    \underbrace{\rest{\ii_T (u_T;J)}{[\pi_{\widehat u_T} (x_1+), \pi_{\widehat u_T} (x_2-)]}}_{\mbox{here } K_T = \reali}
    \subseteq
    \underbrace{\rest{\ii_T (\widehat u_T;J)}{[\pi_{\widehat u_T} (x_1+), \pi_{\widehat u_T} (x_2-)]}}_{\mbox{here } K_T = [x_1,x_2]} \,.
  \end{displaymath}
  Let $u_o \in \L\infty (\reali;J)$ be such that $\ss_T u_o = u_T$.
  Since also
  $\rest{\ii_T (\widehat u_T;J)}{[\pi_{\widehat u_T} (x_1+),
    \pi_{\widehat u_T} (x_2-)]} \neq \emptyset$, choose any
  $\widehat u_o \in \rest{\ii_T (\widehat u_T;J)}{[\pi_{\widehat u_T}
    (x_1+), \pi_{\widehat u_T} (x_2-)]}$. With reference
  to~\eqref{eq:30}, there exists a
  $\overline u_o \in \L\infty (\reali; J)$ such that
  $\rest{\overline{u}_o}{[\pi_{\widehat u_T} (x_1+), \pi_{\widehat
      u_T} (x_2-)]} = \widehat u_o$ and
  $\rest{\ss_T {\overline{u}_o}}{[x_1, x_2]} = \widehat u_T$. Define,
  see Figure~\ref{fig:31},
  \begin{equation}
    \label{eq:31}
    w_o (x) \colonequals
    \left\{
      \begin{array}{r@{\qquad}r@{\,}c@{\,}l}
        u_o (x)
        & x
        & <
        & \pi_{\widehat u_T (x_1+)}
        \\
        \widehat u_o (x)
        & x
        & \in
        & [\pi_{\widehat u_T} (x_1+), \pi_{\widehat u_T}
          (x_2-)]
        \\
        u_o (x)
        & x
        & >
        & \pi_{\widehat u_T (x_2-)}
      \end{array}
    \right.
  \end{equation}
  \begin{figure}[!h]
    \begin{center}
      \begin{tikzpicture}[baseline]
        \draw[thick,-stealth] (-5,0) -- (10,0); %
        \draw[thick,-stealth] (0,-1) -- (0,5); %
        \draw (-0.5,5) node {$t$}; %
        \draw (10,-0.5) node {$x$}; %
        \draw (-0.5,3.75) node {$T$}; %
        \draw (-0.5, -0.5) node {$0$}; %
        \draw[dashed] (-5,4) -- (10,4); %
        \draw (-2.5,0) -- (-1,4); %
        \draw (-3, -0.5) node {$\pi_{\widehat u_T} (x_1+)$}; %
        \draw (-1,4.5) node {$x_1$}; %
        \draw (7,0) -- (5.5,4); %
        \draw (6.5, -0.5) node {$\pi_{\widehat u_T} (x_2-)$}; %
        \draw (5.5,4.5) node {$x_2$}; %
        \draw (2.5,0.25) node {$w_o (x) {=} \widehat u_o (x)$}; %
        \draw (2.25,4.25) node {$\ss_T w_o (x) {=} u_T (x)$}; %
        \draw (8.25, 2) node {$\ss_t w_o (x) {=} \ss_t u_o (x)$}; %
        \draw (-3.35, 2.75) node {$\ss_t w_o (x) {=} \ss_t u_o (x)$};
      \end{tikzpicture}
    \end{center}
    \caption{Notation used in~\eqref{eq:31}.}
    \label{fig:31}
  \end{figure}
  and apply twice Lemma~\ref{lem:glue} using the equalities
  \begin{displaymath}
    \begin{array}{c}
    \ss_T (u_o) (x_1+)
    =
    \widehat u_T (x_1+)
    =
      \ss_T (\overline u_o) (x_1+)
      \\
      \ss_T (u_o) (x_2-)
    =
    \widehat u_T (x_2-)
    =
    \ss_T (\overline u_o) (x_2-)
    \end{array}
  \end{displaymath}
  in~\eqref{eq:24}.  One thus obtains $\ss_T w_o = u_T$, showing that
  $\widehat u_o \in \rest{\ii_T (u_T;J)}{[\pi_{\widehat u_T} (x_1+),
    \pi_{\widehat u_T} (x_2-)]}$.
\end{proofof}

\medskip

\noindent\textbf{Acknowledgment:} The first
author acknowledges the PRIN 2022 project \emph{Modeling, Control and
  Games through Partial Differential Equations} (D53D23005620006),
funded by the European Union - Next Generation EU.  The second author
was supported by the ANR project COSS \emph{COntrol on Stratified
  Structures} (ANR-22-CE40-0010) and from the \emph{Fund for
  International Cooperation} from the University of Brescia.

{

  \small

  \bibliography{xtTraffic}

  \bibliographystyle{abbrv}

}

\end{document}